\documentclass[preprint,12pt]{elsarticle}

\usepackage{amssymb}
\usepackage{lineno}

\usepackage[T1]{fontenc}    
\usepackage{url}            
\usepackage{booktabs}       
\usepackage{amsfonts}       
\usepackage{hyperref}       
\usepackage{nicefrac}       
\usepackage{microtype}      
\usepackage{lipsum}
\usepackage{amsmath}
\usepackage{multicol}
\usepackage{mfirstuc}

\usepackage{lipsum}
\makeatletter
 \def\ps@pprintTitle{%
 \let\@oddhead\@empty
 \let\@evenhead\@empty
 \def\@oddfoot{}%
 \let\@evenfoot\@oddfoot}
\makeatother

\begin{document}

\begin{frontmatter}


\title{Towards Proving the Twin Prime Conjecture using a Novel Method For Finding the Next Prime Number ${P_{N+1}}$ after a Given Prime Number ${P_{N}}$ and a Refinement on the Maximal Bounded Prime Gap ${G_{i}}$}

\author{
  Reema Joshi\\
  Department of Computer Science and Engineering\\
  Tezpur University\\
  Tezpur, Assam 784028, India  \\
  \texttt{reema77@tezu.ernet.in} 
}

\begin{abstract}
This paper introduces a new method to find the next prime number after a given prime ${P}$. The proposed method is used to derive a system of inequalities, that serve as constraints which should be satisfied by all primes whose successor is a twin prime. Twin primes are primes having a prime gap of ${2}$. The pairs ${(5,7),(11,13),(41,43)}$, etcetera are all twin primes. This paper envisions that if the proposed system of inequalities can be proven to have infinite solutions, the Twin Prime Conjecture will evidently be proven true. The paper also derives a novel upper bound on the prime gap, ${G_{i}}$ between ${P_{i} \; and \; P_{i+1}}$, as a function of ${P_{i}}$.
\end{abstract}

\end{frontmatter}


\section{Introduction}
A positive integer ${P>1}$ is called a prime number if its only positive divisors are ${1}$ and ${P}$ itself. For example, 2,3,7,11,19 and so on. There are infinitely many primes. A simple yet elegant proof of this proposition was given by Euclid around 300 B.C. 
One interesting question would be, "does there exist a method for generating successive prime numbers along the infinite pool of primes?". The answer is yes, and is elaborated in
a very well known method for the same, known as the \textbf{Sieve of Eratosthenes}, an algorithm for generating prime numbers up to a given upper bound ${N}$. The complexity 
 of the Sieve of Eratosthenes is ${O(N\log{\log{N}})}$. There are several variants of the Sieve method, which reduce the asymptotic complexity down to ${O(N)}$, like the \emph{Sieve of Sundaram}\cite{aiyar1934sundaram}, \emph{Sieve of Atkin}\cite{atkin2004prime}, and different \emph{Wheel Sieve} methods.\\
This paper proposes a simple method which can accurately calculate the next prime number after a given prime. It may be mentioned beforehand that this method is not very time-efficient(${O(N^2)}$) as compared to existing methods, though its space complexity is ${O(N)}$. This does not make it favourable for prime number generation. However, the goal of this paper is to use the method to revise the upper bound on the prime gap function, ${G_{i}}$, and analyze significant properties of twin primes. As mentioned earlier, the condition for ${\left(P_{i},P_{i+1}\right)}$ to be twin primes is  
\begin{equation*}
    P_{i+1}=P_{i}+2.
\end{equation*}
The \textbf{Twin Prime Conjecture} asserts the infinitude of twin primes and has two different versions. This paper elaborates on the first, as mentioned in \cite{guy1994gaps}, and given below:\\
\textbf{Twin Prime Conjecture}:\emph{ There exist infinitely many primes ${P}$ such that ${P+2}$ is a prime}.\\ In other words,
\begin{equation*}
    \lim_{n \to \infty}(P_{i+1}-P_{i})=2.
\end{equation*}
Though mathematicians hypothesize that the conjecture may well be true\cite{shanks2001solved}, it has not been successfully proven yet.
The distribution of primes seems arbitrary, in some cases, intervals of ${2}$ are seen, for example, the primes ${5}$ and ${7}$ have a gap of 2, while in other cases there are arbitrarily long gaps between consecutive primes\cite{nazardonyavi2012some}. The \textbf{Prime Number Theorem} places a bound on the prime counting function. The \textbf{prime counting function}, denoted by ${\pi_{N}}$ counts the number of primes less than or equal to ${N}$.
Gauss in 1792, stated that,
\begin{equation*}
    \pi_{N} \approx \frac{N}{\log {N}}.
\end{equation*}        
This means for a large enough ${N}$, the probability that any number not exceeding N is prime is close to ${(\frac{1}{\log {N}})}$.
The theorem was independently proven by Jacques Hadamard and  Charles Jean de la Vallée Poussin in 1896.
\\
A detailed mention regarding the bounds on twin primes can be found in \cite{klyve2007explicit}. After the conjecture gained attention, several Sieve based methods have been used at attempts to prove it\cite{ribenboim2012new}.
Zhang, in 2014, made a significant breakthrough along the lines of the Twin Prime Conjecture\cite{zhang2014bounded}. His paper stated that,
\begin{equation*}
    \lim_{n \to \infty}(P_{n+1}-P_{n})<7\times10^{7}.
    \tag{1}
\end{equation*}
This implies that there are an infinite number of primes, ${P_n}$ such that the gap between two consecutive primes is less than 70 million.
In this paper, however, a different approach has been used towards a proof of the Twin Prime Conjecture, based on the proposed method for calculating the successor of a given prime, ${P}$.
The goals of this paper are summarized below:
\begin{enumerate}
    \item To present a novel method for finding ${P_{i+1} \; \textup{given} \; P_{i}}$.
    \item To provide a definite upper bound on the prime gap, ${G_{i}}$ as a function of the ${i^{th} \; \textup{prime}, \;P_{i}}$.
    \item To use the proposed method for setting constraints that should be mandatorily fulfilled by all prime numbers ${P_{i}}$ whose successor prime is its twin, that is, ${P_{i+1}=P_{i}+2}$.
\end{enumerate}

\section{Bounds on the Prime Gap, \texorpdfstring{${G_{i}}$}{}} A prime gap is the difference between two successive primes. If the length of the gap is ${n}$, it implies that there is a sequence of ${n-1}$ successive composite numbers between ${P_{i} \; \textup{and} \; P_{i+1}}$. So,
\begin{equation*}
    G_{i}=P_{i+1}-P_{i}=n,
\end{equation*}
where ${G_{i}}$ represents the prime gap.
${n}$ is always even, except in the case of ${(2,3)}$ where ${n=1}$.\\
We can easily infer that ${G_{1}=1,\;G_{2}=2\;,G_{3}=2}$,${\;G_{4}=11}$ \textup{and so on}. Also, the sequence, ${S}$, 
\begin{equation*}
   S= n!+2, n!+3,n!+4,...,n!+n
\end{equation*}
is a run of ${n-1}$ consecutive composite integers, ${n!+x}$ (${n!+x}$ being divisible by ${x}$). The above equation is evidently part of some prime gap, ${G_{i}}$, such that ${G_{i} \geq n}$. As stated by the Prime Number Theorem, the average prime gap increases with the natural logarithm of the prime number, ${P}$. The ratio ${\dfrac{G_{i}}{\log P_{i}}}$ is called the merit of the gap ${G_{i}}$. The largest known prime gap is 1113106 digits long, and has merit of 25.90. It was discovered by J. K. Andersen,  M. Jansen and P. Cami. Another important term is the \emph{maximal gap}. For two prime gaps ${G_{i} \; \textup{and} \; G_{j}  \;, i,j \in [1,\infty)}$, ${G_{j}}$ is a maximal gap if ${G_{i}<G_{j}\; \forall\; i<j}$. Using exhaustive search, the largest known maximal prime gap having length ${1550}$ was found by Bertil Nyman\cite{nyman2003new} and occurs after the prime number ${18361375334787046697}$.\\
\textbf{Bertrand's Postulate} states that for every integer ${k\left( k >1\right)}$, there exists a prime number ${P}$, such that,
\begin{equation*}
    k < P < 2k.
\end{equation*}
Another form of the postulate is that, for the ${i^{th}}$ prime, ${P_{i}}$,
\begin{equation}
    P_{i+1} < 2P_{i}.\tag{2}
\end{equation}
This implies that ${G_{i} < P_{i}}$.
It was proposed as a conjecture initially, and proved by  Chebyshev. This postulate is definitely weaker than the Prime Number Theorem, since there are roughly ${\dfrac{k}{\ln k}}$ primes upto ${k}$. This count is greater than that stated by Bertrand's postulate. An upper bound can be inferred on the gaps between primes using the Prime Number Theorem, as follows:\\
\emph{For every ${\epsilon>0}$, there is a number ${N}$ such that ${\forall \;n>N,\; G_{i}<P_{i}\epsilon}$}.
Also,
\begin{equation*}
    \lim_{n \to \infty}\dfrac{G_{i}}{P_{i}}=0
\end{equation*}
because the gaps get arbitrarily small compared to the prime numbers.
On the other hand, \textbf{Legendre's Conjecture} states that for every integer ${k, k >1}$, there exists a prime number ${P}$, such that,
\begin{equation*}
    k^2 < P < \left(k+1\right)^2.
\end{equation*}
\textbf{Cramer's Conjecture} This conjecture, devised by Harald Cramer in 1936\cite{cramer1936order} provides an asymptotic bound on prime gaps:
\begin{equation*}
    P_{i+1}-P_{i}=O\left(\left(\log P_{i}\right)^2\right) 
    \tag{3}
\end{equation*}
where ${O}$ represents the Big-Oh notation. There are several versions of Cramer's Conjecture, of which (3) has not been proven yet. There are also weaker versions of this conjecture, which have been verified using conditional proofs, listed below:
\begin{enumerate}
    \item ${P_{i+1}-P_{i}=O\left(\sqrt{P_{i}} \log P_{i}\right)}$, which Cramer proved assuming the Reimann Hypothesis\cite{cramer1936order}.
    \item ${P_{i+1}-P_{i}=O\left(P_{i}^{0.525}\right)}$ proved by Baker, Harman and Pintz\cite{baker1996difference}.
    \item With reference to the previous two points, prime gaps can grow more than logarithmically, therefore, the upper bound would have to be further worked on. It was proved that 
    \begin{equation*}
        \limsup_{n \to \infty}\dfrac{P_{i+1}-P_{i}}{\log P_{i}}=\infty.
    \end{equation*}
    \item Point 3 was improved as under:
    \begin{equation*}
        \limsup_{n \to \infty}\dfrac{P_{i+1}-P_{i}}{\log P_{i}}\cdot\dfrac{\left(\log \log \log P_{i}\right)^2}{\log \log P_{i} \log \log \log \log P_{i}}>0.
    \end{equation*}
\end{enumerate}
It was conjectured by Paul Erdos that the left hand side of the above equation is infinite, and this was successfully proven in \cite{ford2016large}.\\
A strong form of Cramer's conjecture was provided by Daniel Shanks in \cite{shanks1964maximal} as 
\begin{equation*}
    G_{i} \approx \log ^2 P_{i}
\end{equation*}
for the ${i^{th}}$ prime ${P_{i}}$.
The maximal gap ${G_{i}}$ has been expressed in terms of the prime counting function as elaborated by Wolf in \cite{wolf2014nearest},
\begin{equation*}
    G_{i} \approx \dfrac{P_{i}}{\pi\left(P_{i}\right)}\left(2 \log \pi\left(P_{i}\right) - \log \left(P_{i}\right) + c\right)
\end{equation*}
where ${c= \log {C_{2}}=0.2778769... \;\textup{and} \; C_{2}=1.3203236}$ is double the twin primes constant.
Gauss's approximation stated earlier, 
\begin{equation*}
    \pi_{N} \approx \frac{N}{\log {N}},
\end{equation*}   
would reduce the previous equation to
\begin{equation*}
    G_{i} \approx \log \left(P_{i}\right)\left(\log \left(P_{i}\right) - 2\log \log \left(P_{i}\right)\right),
\end{equation*}
which for large values of ${P_{i}}$ reduces to Cramer's Conjecture, ${G_{i} \approx \log^2P_{i}}$.
Zhang\cite{zhang2014bounded} proved that there are an infinite number of primes which have a prime gap less than 70,000(Eq. 1). After this, James Maynard reduced the bound \cite{maynard2015small} by proving that 
\begin{equation*}
       \lim_{n \to \infty}(P_{n+1}-P_{n})\leq600.
\end{equation*}
The Polymath Project brought mathematicians to work together on extending Zhang's proof, and finally, the bound has been improvised and brought down to 246.
The following sections discuss the working of the proposed method, its direct conclusions on the bounds of the prime gap, ${G_{i}}$ and how it paves way for an alternative approach towards validating the truth of the Twin Prime Conjecture.

\section{Method to Find \texorpdfstring{${P_{N+1}}$}{} Using \texorpdfstring{${P_{N}}$}{}}
The proposed method revolves around a concept named here as \emph{“Slack”}. The idea of a slack can be understood as follows:\\
In any division ${\dfrac{P}{D}}$, where ${P}$ and ${D}$ are both integers, there are two cases that can occur: Either ${P}$ is divisible by D, or it is not. If not, we need to calculate a value ${S}$, such that after adding ${S}$ to ${P}$, we get ${N=P+S}$, where ${N}$ is divisible by D. So, calculation of the slack is the process finding the amount by which ${P}$ lacks in order to be divisible by ${D}$.
For example, 8 is not divisible by 3. However, adding 1 to 8 gives 9, which clearly is divisible
by 3. So, in this case, the slack, S= 1.
Can we devise a formula for the slack? Yes.
The slack, ${S}$ is calculated using:
\begin{equation}
  S=D\cdot\left(\lfloor{\dfrac{P}{D}}\rfloor+1\right)-P. \tag{4}
\end{equation}
Taking an example to illustrate this concept, let ${P = 7 \;, D=3}$. Then,
\begin{align*}
S&=D\cdot\left(\lfloor{\dfrac{P}{D}}\rfloor+1\right)-P \\
 &= 3\cdot\left(\lfloor{\dfrac{7}{3}}\rfloor+1\right)-7\\
 &= 3\cdot(2+1)-7\\
 &= 9-7\\
 &= 2.
\end{align*}
 Adding the obtained slack, ${S=2}$ to the dividend 7 gives ${N=7+2=9}$, which is divisible by 3.
\section{Algorithm to Find the Next Prime Number}
This section explains the proposed algorithm in detail.
\\
\textbf{Goal}: To find the successor of the ${i^{th}}$ prime number, ${P_{i}}$, where ${i\;\in[1,\infty) }$. 
The first three primes are statically initialized to ${P_{1} = 2, P_{2} = 3}$ and ${P_{3}=5}$(the method will only generate primes ${>5)}$.
The following are the steps to generate ${P_{i+1}}$ given ${P_{i}}$:
\begin{enumerate}
    \item Start with a given prime number, assign it to ${P_{i}}$.
    \item Find all potential divisors of ${P_{i}}$ , starting from 2 up to ${\dfrac{P_{i}-1}{2}}$ and store this in a list, say, ${divisorList}$ for ${P_{i}}$.
    \item For each ${D}$ in ${divisorList}$,
    \begin{enumerate}
        \item If ${P_i}$ is divisible by ${D}$, do nothing,\\
        else
        \item Store the slack S for the combination ${(P_i,D)}$ into a list of slacks, say, ${slackList}$ for ${P_{i}}$.
    \end{enumerate}
    \item Go through ${slackList}$, and identify the first even number ${E}$, missing in the list.
    \begin{enumerate}
        \item The value of the slack falls in the range ${[1,{\dfrac{P_{i}-1}{2}-1}]}$, since ${D}$ lies in the range ${[2,\dfrac{P_{i}-1}{2}]}$, and the slack will never be ${0}$ since prime numbers will not be divisible by any ${D}$.
        \item The even numbers in the slack list vary from ${2}$ to ${\dfrac{P_{i}-1}{2}-1}$, if ${\dfrac{P_{i}-1}{2}}$ is odd, else from 2 to ${\dfrac{P_{i}-1}{2}-2}$, if ${\dfrac{P_{i}-1}{2}}$ is even.
        \item In case no even slacks in the range given by point b are missing in ${slackList}$, then we go beyond the range to pick the next even number, ${E=\dfrac{P_{i}-1}{2}}$, if ${\dfrac{P_{i}-1}{2}}$ is even, else ${E=\dfrac{P_{i}-1}{2}+1}$, if ${\dfrac{P_{i}-1}{2}}$ is odd.
        \item Add ${E}$ to ${P_i}$ to obtain the next prime number, ${P_{i+1}}$.
    \end{enumerate}
\end{enumerate}
The main idea behind this algorithm is to find a slack which when added to the current prime number ${P_i}$ will still not make it divisible by any of the potential divisors, therefore, this number may also be a prime number. However, since ${P_{i+1}}$ will have more divisors than ${P_{i}}$(including divisors of ${P_{i}}$), it is required to be sure that ${P_{i+1}}$ is not divisible by them all. In the next section, a mathematical proof of the proposed method is provided to show that the generated prime, ${P_{i+1}}$,  will not be divisible by any of its potential divisors. Also, the infinitude of primes signifies that there will always be a next prime number ${P_{i+1}}$ after any ${i^{th}}$ prime ${P_{i}}$. The following examples cover the different cases discussed in point 4 of the proposed algorithm.\\
\textbf{\emph{Example 1}}.  Let ${P=11}$. For prime number ${P_{i}=11},\;i=5\;$since ${11}$ is the ${5^{th}}$ prime number. The potential divisors of ${11}$ are ${2}$ to  ${5\left(\dfrac{11-1}{2}=5\right)}$. The values of the slacks corresponding to each divisor are given in Table 1.
\begin{table}[ht]
\centering
\begin{tabular}[t]{c c }
\hline
    Divisor & Slack \\
\hline    
    2 & 1 \\
\hline    
    3 & 1 \\
\hline    
    4 & 1 \\
\hline    
    5 & 4 \\
\hline    
\end{tabular}
\caption{Slack List For ${P_{i}=11}$}
\label{tab:Table 1}
\end{table}
Now, from the list of slacks, the first missing even number is ${2}$, so ${E=2}$. Therefore, next prime number ${P_{i+1}= P_{i} + E =11 + 2 = 13}$. 
Verifying point 4(b), ${E}$ lies in the interval ${[2,\dfrac{11-1}{2}-1]}$, that is ${[2,4]}$, since ${\dfrac{P_{i}-1}{2}}$ is odd.\\
\textbf{\emph{Example 2}}.  Let ${P=29}$. For prime number ${P_{i}=29}$, the potential divisors are 2 to  ${14\left(\dfrac{29-1}{2}=14\right)}$. 
\begin{table}[ht]
\centering
\begin{tabular}[t]{c c}
\hline
    Divisor & Slack \\
\hline    
    2 & 1 \\
\hline    
    3 & 1 \\
\hline    
    4 & 3 \\
\hline    
    5 & 1 \\
\hline  
    6 & 1 \\
\hline
    7 & 6 \\
\hline 
    8 & 3 \\
\hline  
    9 & 7 \\
\hline  
    10 & 1 \\
\hline
    11 & 4 \\
\hline  
    12 & 7 \\
\hline  
    13 & 10 \\
\hline  
    14 & 13 \\
\hline  
\end{tabular}\\
\caption{Slack List For ${P_{i}=29}$}
\label{tab:Table 2}
\end{table}
Looking at Table 2, from the list of slacks, the first missing even number is ${2}$, so ${E=2}$.
Therefore, next prime number ${P_{i+1}= P_{i} + E =29 + 2 = 31}$. 
Verifying point 4(b), ${E}$ lies in the interval ${[2,\dfrac{29-1}{2}-2]}$ since ${\dfrac{P_{i}-1}{2}}$ is even, that is, ${[2,12]}$.
\\\textbf{\emph{Example 3}}. Let ${P=7}$. For prime number ${P_{i}=7}$, the potential divisors are ${2}$ to  ${3\left(\dfrac{7-1}{2}=3\right)}$. 
\begin{table}[ht]
\centering
\begin{tabular}[t]{c c}
\hline
    Divisor & Slack \\
\hline    
    2 & 1 \\
\hline    
    3 & 2 \\
\hline    
\end{tabular}
\caption{Slack List For ${P_{i}=7}$}
\label{tab:Table 3}
\end{table}
From Table 3, we find that the slack list contains ${1,2}$. Looking at point 4(b) of the algorithm, ${E}$ should be lying in the interval ${[2,\dfrac{7-1}{2}-1]}$ since ${\dfrac{P_{i}-1}{2}}$ is odd, that is, ${[2,2]}$.
But ${2}$ exists in the slack list, and we cannot choose it.
So, we go beyond this range to pick the next even number after ${2}$, getting ${4}$.
Thus, ${E=4}$ and next prime number ${P_{i+1}=P_{i}+E=7+4=11}$.
This verifies point 4(c), because ${\dfrac{P_{i}-1}{2}=\dfrac{7-1}{2}=3}$  is odd, 
${E=\dfrac{P_{i}-1}{2}+1}$, so ${E=3+1=4}$.
\\\textbf{\emph{Example 4}}. Let ${P=5}$. For prime number ${P_{i}=5}$, the potential divisors are ${2}$ to  ${\dfrac{5-1}{2}=2}$. There is only one divisor, ${2}$ in this case. 
\begin{table}[ht]
\centering
\begin{tabular}[t]{c c}
\hline
    Divisor & Slack \\
\hline    
    2 & 1 \\
\hline    
\end{tabular}
\caption{Slack List For ${P_{i}=5}$}  
\label{tab:Table 4}
\end{table}
Table 4 does not contain any even slack because the first even number, ${2}$ is the only divisor, and so the maximum and only slack is 1.
Since the minimum value of ${E}$ for any prime number has to be 2(${E}$ is even), and 2 will obviously not exist in the slack list, we choose the first missing even number beyond the slack list, to get ${E=2}$.
Next prime number, ${P_{i+1} = P_{i}+E=5+2=7}$. This verifies point 4(c).

\subsection{Verification of the Proposed Method}
The proposed method has been tested with the help of a computer program and successfully verified for the first 1,00,000 primes. Further, in the next section, a mathematical proof is provided to show that the method will reliably work for all successive prime numbers.

\subsection{Proof of the Proposed method}
According to the \textbf{Prime Number Theorem}, on an average, the gap between a given prime number ${N}$ and its successor prime number increases along ${\log N}$\cite{klyve2007explicit}.
Let ${P_{i}}$ be the ${i^{th}}$ prime, where ${i\in[1,\infty)}$ and ${P_{i+1}}$ be ${P_{i}}$'s successor prime.
Since the potential divisors of ${P_{i}}$ vary from ${2}$ to  ${\dfrac{P_{i}-1}{2}}$, the range of the slacks varies from ${1}$ to ${\dfrac{P_{i}-1}{2}-1}$ . 
Since only an even slack is added to ${P_{i}}$ to generate the next prime number ${P_{i+1}}$, it can be said that the range of ${P_{i+1}}$ varies from ${P_{i}+2}$  to  ${P_{i}+ \left(\dfrac{P_{i}-1}{2}-1\right)}$, if ${\dfrac{P_{i}-1}{2}-1}$ is even, else ${P_{i}+2}$  to ${P_{i}+ \left(\dfrac{P_{i}-1}{2}-2\right)}$, if ${\dfrac{P_{i}-1}{2}-1}$ is odd.
Thus, minimum value of ${P_{i+1}= P_{i}+ 2}$.
Also, we need to consider two more cases from point 4(c), which says, ”In case no even slacks in the range given by point 4(b) are missing, the we go beyond the range to pick the next even number, that is, ${E=\dfrac{P_{i}-1}{2}}$ if ${\dfrac{P_{i}-1}{2}}$ is even, else ${E=\dfrac{P_{i}-1}{2}+1}$ if ${\dfrac{P_{i}-1}{2}}$ is odd”. So, two other possible ranges for ${P_{i+1}}$ are ${P_{i}+2}$  to ${P_{i}+\dfrac{P_{i}-1}{2}+1}$ and ${P_{i}+2}$  to ${P_{i}+\dfrac{P_{i}-1}{2}}$.
The maximum value of ${P_{i+1}}$, keeping in mind point 4(c) where slack exceeds normal range, is
\begin{equation}
    P_{i+1} = P_{i} + \dfrac{P_{i}-1}{2}+1
            = \dfrac{3P_{i}}{2}+\dfrac{1}{2}. \tag{5}
\end{equation}
As discussed earlier, for ${P_{i+1}}$ (obtained using this method) to be part of a twin prime pair, ${P_{i+1}}$ should not be divisible by any of its potential divisors.
These will be  2, 3, ...,  ${\dfrac{P_{i}-1}{2}, \dfrac{P_{i}-1}{2}+1,..., \dfrac{P_{i+1}-1}{2}}$.
Since the maximum value of ${P_{i+1}}$ is ${P_{i}+\left(\dfrac{P_{i}-1}{2}+1\right)}$, the maximum value of ${P_{i+1}}$'s divisor is, ${\dfrac{P_{i+1}-1}{2}}$ which is equivalent to
\begin{equation}
    \dfrac{P_{i}+\dfrac{P_{i}-1}{2}+1-1}{2}
    =\dfrac{3P_{i}}{4}-\dfrac{1}{4}. \tag{6}
\end{equation}
The maximum divisor for ${P_{i+1}}$ may be less than ${\dfrac{3P_{i}}{4}-\dfrac{1}{4}}$, in cases where ${P_{i+1}<\left(P_{i}+\dfrac{P_{i}-1}{2}+1\right)}$. 
For example, If ${P_{i}=7}$ and ${P_{i+1}=11}$, slack ${E=4}$. Now, maximum divisor of ${P_{i+1}}$, that is, ${11=\dfrac{3\times7}{4}-\dfrac{1}{4}=5}$. It satisfies relation (3) because the slack ${4}$ equals  ${\dfrac{P_{i}-1}{2}+1}$. \\
But if ${P_{i}=11}$ and ${P_{i+1}=13}$, maximum divisor of ${P_{i+1}}$, that is, ${\dfrac{3P_{i}}{4}-\dfrac{1}{4}=8}$. Of course, ${8}$ cannot be a potential divisor of ${13}$, this is happening because we are considering the maximum possible value of ${P_{i+1}}$ considering cases of all prime number successors. In this case relation (3) is not satisfied, since the slack ${2}$ $\neq$ ${\dfrac{P_{i}-1}{2}+1}$.
This implies that it depends where, out of the four possible ranges mentioned in points 4(b) and 4(c) of the algorithm does the value of the slack fall. This is not a problem because as discussed, we have considered the maximum possible divisor value out of all the four cases. In any case, the maximum divisor for ${P_{i+1}=\dfrac{P_{i+1}-1}{2}}$.
Now, the divisors from ${2}$ up to   ${\dfrac{P_{i}-1}{2}}$   will not be able to divide ${P_{i+1}}$ because ${P_{i+1}}$ is generated from the missing slack, that is, the slack value which when added to ${P_{i}}$ will still not make the result divisible by the potential divisors of ${P_{i}}$(which range from ${2}$ to ${\dfrac{P_{i}-1}{2}}$). But what about the remaining divisors of ${P_{i+1}}$ which are beyond the divisor list of ${P_{i}}$ (since ${P_{i+1}}$>${P_{i}}$ it will have more divisors than ${P_{i}}$)? The next section proves that even these remaining divisors cannot divide ${P_{i+1}}$. 
\subsubsection{Mathematical Proof}
Remaining divisors of ${P_{i+1}}$(apart from common divisors  of ${(P_{i},P_{i+1}))}$  will be : ${\dfrac{P_{i}-1}{2}+1}$ up to  ${\dfrac{3P_{i}}{4}-\dfrac{1}{4}}$ (maximum divisor upper limit for ${P_{i+1}}$).\\
Let the symbol ${x}$ represent any of these remaining divisors of ${P_{i+1}}$. The following observations are important:
\begin{enumerate}
    \item Since ${P_{i+1}}$ is an odd number, 
    \begin{equation*}
        P_{i+1} \neq 2x,
        \forall x \; where\\ \; x \in [\dfrac{P_{i}-1}{2}+1,\dfrac{3P_{i}}{4}-\dfrac{1}{4}].
    \end{equation*}
    \item Now let us take the smallest of these remaining divisors,which is the minimum value of ${x}$ from its range, that is, ${x =\dfrac{P_{i}-1}{2}+1}$. Assuming ${P_{i+1} = 3x}$, we get
\begin{equation}
    P_{i+1} = 3\cdot\left(\dfrac{P_{i}-1}{2}+1\right)
            = \dfrac{3P_{i}}{2}-\left(\dfrac{3}{2}+3\right)
            = \dfrac{3P_{i}}{2}+\frac{3}{2}.\tag{7}
\end{equation}
\end{enumerate}
This violates relation (5) because maximum value of ${P_{i+1}}$ should be ${\dfrac{3P_{i}}{4}+\dfrac{1}{2}}$, which is clearly less than that evaluated by (7). So,
(7) contradicts (5). Since a contradiction is encountered, our initial assumption ${P_{i+1} = 3x}$ is wrong.
As the minimum divisor ${x=\dfrac{P_{i}-1}{2}+1}$ had its lowest multiple(${3x}$) exceed ${P_{i+1}}$, it is evident that ${P_{i+1}}$ will not be divisible by any of the successive divisors.
Also, we have considered the maximum value of ${P_{i+1}}$ which is ${P_{i}+\left(\dfrac{P_{i}-1}{2}+1\right)}$. This itself is below ${3x}$, where ${x}$ is the lowest value of the divisor, ${x=\dfrac{P_{i}-1}{2}+1}$.
The following is a list all possible values that ${P_{i+1}}$ can take: 
\begin{enumerate}
    \item 	${P_{i+1}=P_{i}+\left(\dfrac{P_{i}-1}{2}+1\right)}$.
    \item   ${P_{i+1}=P_{i}+\dfrac{P_{i}-1}{2}}$.
    \item  	${P_{i+1}=P_{i}+\left(\dfrac{P_{i}-1}{2}-1\right)}$.
    \item  	${P_{i+1}=P_{i}+\left(\dfrac{P_{i}-1}{2}-2\right)}$, and
    \item   All other values of ${P_{i+1}}$ lower than given in the above four points, when the slacks are below ${\dfrac{P_{i}-1}{2}}$,
    ${\dfrac{P_{i}-1}{2}-1}$, ${\dfrac{P_{i}-1}{2}+1}$, ${\dfrac{P_{i}-1}{2}-2}$, as follows:
    \begin{enumerate}
       \item   ${P_{i+1} < P_{i}+\left(\dfrac{P_{i}-1}{2}+1\right)}$.
       \item   ${P_{i+1} < P_{i}+\dfrac{P_{i}-1}{2}}$.
       \item   ${P_{i+1} < P_{i}+\left(\dfrac{P_{i}-1}{2}-1\right)}$.
       \item   ${P_{i+1} < P_{i}+\left(\dfrac{P_{i}-1}{2}-2\right)}$.
    \end{enumerate}
\end{enumerate}
It can be understood that the other possible values of ${P_{i+1}}$ listed in points 2, 3, 4 and 5 are lower than the value of ${P_{i+1}}$ given in point 1. So they will all be less than ${3x}$. 
Therefore, none of the remaining divisors from ${\dfrac{P_{i}-1}{2}+1}$ up to ${\dfrac{P_{i+1}-1}{2}}$ will be able to divide ${P_{i+1}}$.
Thus, ${P_{i+1}}$ will always turn out to be a prime number.

\subsection{Results on the Maximum Prime Gap Between \texorpdfstring{{}${P_{i} \; \textup{and} \; P_{i+1}}$}{}}
From the discussion in point 4(b) of the algorithm, the minimum value of ${E=2}$, where ${E}$ is the missing slack, added to ${P_{i}}$ to obtain ${P_{i+1}}$. For example, the case of ${(3,5)}$. What about the maximum value of ${E}$? The answer would give us information about the upper bound on the prime gap function, ${G_{i}}$. Going by the Prime Number Theorem, average prime gap between any prime ${N}$ and its successor prime is ${\log N}$, while Bertrand's postulate bounds the maximum prime gap to ${N-1}$(eq. 2). This paper further shrinks the upper bound on ${G_{i}}$ to ${\dfrac{P_{i}+1}{2}}$. The proof is explained below:
\paragraph{Proof}From eq.(5), we infer that the maximum value of ${P_{i+1}}$ is
\begin{equation*}
    P_{i+1}=P_{i}+\left(\dfrac{P_{i}-1}{2}+1\right).
\end{equation*}
This implies that the maximum difference between ${P_{i}}$ and ${P_{i+1}}$ is
\begin{equation*}
    G_{i} =P_{i+1}-P_{i} = \left(\dfrac{P_{i}-1}{2}+1\right)
= \dfrac{P_{i}}{2}+\dfrac{1}{2}
= \dfrac{P_{i}+1}{2}.
\end{equation*}
${P_{i}+1}$ will always be an even number, so ${\dfrac{P_{i}+1}{2}}$ will be an integer, except in the special case where ${P_{i}=2}$, we find the prime gap between ${P_{i} \; and \; P_{i+1}, \textup{that is}}$, ${(2,3)}$ equals ${1.5}$. Therefore, it is important to round this expression to an integer. The floor function will be used for rounding since the prime gap between ${2}$ and ${3}$ equals ${1}$, and for all other primes,${\dfrac{P_{i}+1}{2} \in \mathbb{Z}^+}$ and is ${>1}$. Therefore
\begin{equation*}
    1 \leq G_{i} \leq \lfloor\dfrac{P_{i}+1}{2}\rfloor.
    \tag{8}
\end{equation*}
The prime gap between ${2 \; \textup{and} \; 3}$, ${G_{1}}$, satisfies (8), as
\begin{align*}
  G_{1}  &=\lfloor\left(\dfrac{2+1}{2}\right)\rfloor\\
    &=1.
\end{align*}
Eq.(8) can be verified for any prime pair, ${\left(P_{i},P_{i+1}\right)}$, also because the proposed method has been mathematically proved in Section 4.2.1. Further, (8) also fulfils Bertrand's postulate since ${\lfloor\dfrac{P_{i}+1}{2}\rfloor < P_{i}}$, in fact, the upper  bound on ${G_{i}}$ is reduced by approximately ${\dfrac{1}{2}}$ as compared to the postulate.

\section{Towards Proving the Twin Prime Conjecture}
According to (4), the formula to calculate the slack in the case of any division ${\dfrac{P}{D}}$ is
\begin{align*}
    Slack &= D\cdot\left(\lfloor{\dfrac{P}{D}}\rfloor+1\right)-P\\
          &= D\cdot\left(\frac{P}{D}-\alpha+1\right)-P, \tag{9}
\end{align*}
\textup{where} ${\alpha  \, \in [0,1)}$.\\
We now analyse what happens in the case of the twin prime pair ${\left(P_{i},P_{i+1}\right)}$ .
${P_{i}}$'s potential divisors lie in ${[2, \dfrac{P_{i}-1}{2}]}$. Every division, ${\dfrac{P_{i}}{D_{k}}}$ results in a corresponding slack ${S_{k}}$(${k \in [2,\dfrac{P_{i}-1}{2}]}$) which can take values from ${1}$ up to ${\dfrac{P_{i}-1}{2}-1}$.
Now, for a prime ${P_{i}}$ having a twin successor, ${P_{i+1}}$(${P_{i+1} = P_{i}+2}$), 
\begin{align*}
  S_{k} \neq 2 \; \forall \; k \in [2,\dfrac{P_{i}-1}{2}].\tag{10}
\end{align*}
Eq.(10) is required because if the slack is ${2}$, then ${P_{i}+2}$ will not be a prime number, because ${2}$ will not be missing in the slack list(the algorithm says that only the missing slack is added to ${P_{i}}$ to obtain the next prime). In this case, ${P_{i}}$ will not have a twin prime.
So, the following system of inequalities constrain the slack to ensure that it is does not equal 2. All of the following constraints should be simultaneously fulfilled by the first number in any twin prime pair.\\
\textbf{Constraints}:
\begin{equation*}
S_{2}=2\cdot\left(\lfloor{\dfrac{P}{2}}\rfloor+1-P\right)=1.
\tag{${C_{1}}$}
\end{equation*}
${C_{1}}$ will be satisfied by all prime numbers, since all primes(except 2) are odd. Next,
\begin{multicols}{2}
\noindent
    \begin{equation*}
      S_{3}=3\cdot\left(\lfloor{\dfrac{P}{3}}\rfloor+1-P\right)>2,
        \; \textup{or}
    \end{equation*}
     \begin{equation*}
        S_{4}=4\cdot\left(\lfloor{\dfrac{P}{4}}\rfloor+1-P\right)>2, \; \textup{or}
    \end{equation*}
     \begin{equation*}
      \textup{ Moving successively},
        S_{k}=k\cdot\left(\lfloor{\dfrac{P}{k}}\rfloor+1-P\right)>2, \; \textup{or}
    \end{equation*}    
    \begin{equation*}
        S_{3}=3\cdot\left(\lfloor{\dfrac{P}{3}}\rfloor+1-P\right)=1,
    \tag{$C_{2}$}
    \end{equation*}
    \begin{equation*}
        S_{4}=4\cdot\left(\lfloor{\dfrac{P}{4}}\rfloor+1-P\right)=1,
    \tag{$C_{3}$}
    \end{equation*}
    \begin{equation*}
       S_{k}=k\cdot\left(\lfloor{\dfrac{P}{k}}\rfloor+1-P\right)=1,
     \tag{${C_{k-1}}$}
    \end{equation*}
\end{multicols}
where ${k = \dfrac{P_{i}-1}{2}}$.\\   
If ${P_{i}}$ and ${P_{i+1}}$ are twin primes, ${P_{i}}$ will satisfy a particular set of inequalities specific to it, that is, values in  \emph{slackList} of ${P_{i}}$ will either be ${1}$ or ${>2}$. As the algorithm states in \emph{Section 4}, the next prime number ${P_{i+1}}$ will then be calculated by adding ${2}$(the first missing even number in the slack list) to ${P_{i}}$, as,
\begin{equation*}
    P_{i+1}  = P_{i}+2,
\end{equation*}
resulting in the twin prime pair, ${(P_{i},P_{i+1})}$.
The goal of this paper is to present the above set of constraints as an important property of twin primes, and therefore, a requirement that must be followed by all primes whose successor is a twin prime.\\
According to the \textbf{Twin Prime Conjecture}, there are infinite pairs of prime numbers, ${(P_{i},P_{i+1})}$ that differ by ${2}$. To prove this conjecture, one needs to prove that the system of inequalities ${C_1}$ to ${C_{k-1}}$ will have infinite solutions. Now, ${C_1}$ will hold for all prime numbers except ${2}$, since they are all odd, so the slack ${S}$ on dividing any prime number ${P}$ by ${2}$ will always be ${1}$.

\subsection{Alternative System of Inequalities For Twin Primes}
The concept of a remainder and slack are very close to each other. A remainder,${R\equiv P ~(\textup{mod} ~D)}$ is a value that is left over after division, say, ${\dfrac{P}{D}}$, when ${P}$ is not divisible by ${D}$. A slack ${S}$, on the other hand, shows by what amount does ${P}$ lack in order to be divisible by ${D}$.
Therefore,
\begin{equation}
    R = D - S. 
    \tag{11}
\end{equation}
where ${R}$ is the remainder and ${S}$ is the slack when any prime number ${P}$ is divided by a divisor ${D}$.\\
Now, according to (9), 
\begin{align*}
    Slack &= D\cdot\left(\dfrac{P}{D}-\alpha+1\right)-P, \\
\end{align*}
where  ${\alpha  \, \in [0,1)}$.  \\    
Here, ${\alpha}$ can also be expressed as ${\dfrac{R}{D}}$, where ${R}$ and ${D}$ are the remainder and divisor respectively. Going by the \textbf{Remainder Theorem}, the possible values of the remainder ${R}$ lie in the range ${[0,(D-1)]}$. In this case however, the remainder ${R}$ will never be ${0}$, as prime numbers will not be divisible by any divisor ${D}$, where ${D}$ lies in the range ${[2,\dfrac{P-1}{2}]}$. Thus, the range of ${R}$ will be ${[1,D-1]}$.
The expression for ${\alpha}$ is 
\begin{align*}
     \alpha &=\dfrac{R}{D}\\
           &= \dfrac{P}{D}-\lfloor\dfrac{P}{D}\rfloor
         \tag{12}\\
     \textup{and}, \; R&\equiv P ~(\textup{mod} ~D).
      \tag{13}
\end{align*}                                     
 It may be mentioned that  ${\dfrac{P}{D}}$ is a decimal value and ${\lfloor{\dfrac{P}{D}}\rfloor}$ is the greatest integer less than ${\dfrac{P}{D}}$.
Since the remainder ${R}$ goes from ${1}$ to ${D-1}$, the range of ${\alpha}$  is ${[\dfrac{1}{D} ,\dfrac{D-1}{D}]}$.\\
Now, the set of constraints ${C_{1}}$ to ${C_{k-1}}$ were  mentioned as necessary conditions for any prime ${P_{i}}$ to have its successor, ${P_{i+1}}$ be a twin prime.
These constraints were based on the slack, ${S}$.
Another way to express the constraints ${C_{1}}$ to ${C_{k-1}}$ is to center them around ${R}$, instead of ${S}$. Since the \emph{slack} and \emph{remainder} are related as shown in (11), we can do this mapping easily. Therefore, in any twin prime pair, ${(P_{i},P_{i+1})}$, where 
${P_{i+1}=P_{i}+2}$, all the inequalities in the  following set of constraints should be fulfilled simultaneously. If any constraint does not hold, then ${(P_{i},P_{i+1})}$ will not be a twin prime pair.\\ 
\textbf{Constraints:}
\begin{align*}
        R_{2} &\neq D_{2} - 2, 
        \tag{14}\\
        R_{3} &\neq D_{3} - 2, 
        \tag{15}\\
        &... \\
        R_{k} &\neq D_{k} - 2, 
\end{align*}
where ${R_{k}}$ is the remainder obtained on dividing ${P_{i}}$ by ${D_{k}}$ , and ${k}$ lies in the range ${[2,\dfrac{P_{i}-1}{2}]}$.
The above set of inequalities containing ${\dfrac{P_{i}-1}{2}}$ constraints may collectively be called "\textbf{${R-constraints}$}". It is obvious that the larger the value of the prime number, the more the number of inequalities in this set.
Since ${R\equiv P ~(\textup{mod} ~D)}$, as derived in (13), the value of ${R}$ depends on both ${P}$ and ${D}$.
\paragraph{Justification for the constraints.} Eq.(11) mentions that
\begin{equation*}
    R=D-S.
\end{equation*}
If ${R=D-2}$, it means the the slack is 
\begin{align*}
    S &= D-(D-2)\\
      &= 2.
\end{align*}
This is not acceptable according to the proposed algorithm. If any slack from the list of slacks ${S_{1} \;to\; S_{k}}$, ${k=\dfrac{P_{i}-1}{2}}$ contains 2, then 2 will not be the first missing even number in the algorithm. Therefore,
\begin{equation*}
    P_{i+1} \neq P_{i}  + 2.
\end{equation*}
This implies that ${(P_{i},P_{i+1})}$ will not be a twin prime pair.
The following are two example cases which verify the truth of the proposed system of constraints:\\
\textbf{\textbf{Case 1}}. Twin prime pair ${\left(P_{i}, P_{i+1}\right) = {\left(11, 13\right)}}$.\\
The potential divisors of ${11}$ are ${2 \;to\; 5 \left(\dfrac{11-1}{2} = 5\right)}$.
 Below, each constraint is enlisted in the way it should hold for each divisor of the first number in the twin prime pair,which is 11.
\begin{align*}
   	D&=2,D-2=0, R\equiv 11~(\textup{mod} ~2), R \neq D-2,\\
	D&=3,D-2=1, R\equiv 11~(\textup{mod} ~3), R \neq D-2,\\
	D&=4,D-2=2, R\equiv 11~(\textup{mod} ~4), R \neq D-2,\\ D&=5,D-2=3, R\equiv 11~(\textup{mod} ~5), R \neq D-2.
\end{align*}
Since all the inequalities are satisfied by the prime number ${11}$, it means ${11}$ has a twin prime as its successor. Also, since ${R \neq D-2 \; \forall \; (R,D)}$, the first missing slack(denoted by ${E}$) is 2. Therefore, ${Next Prime = 11+E = 11+2 = 13}$, which is correct.

\textbf{\textbf{Case 2}}. A non-twin prime pair ${\left(P_{i}, P_{i+1}\right) = {\left(13, 15\right)}}$. The potential divisors of ${13}$ are ${2 \;to\; 6 \left(\dfrac{13-1}{2} = 6\right)}$.
We enlist each constraint as it should hold for each divisor of the first number in the twin prime pair, which is 13.
\begin{align*}
   	&D=2, D-2=0, R\equiv 13~(\textup{mod} ~2), R \neq D-2,\\
	&D=3, D-2=1, R\equiv 13~(\textup{mod} ~3), R = D-2,\\
	&D=4, D-2=2, R\equiv 13~(\textup{mod} ~4), R \neq D-2, \\ &D=5, D-2=3, R\equiv 13~(\textup{mod} ~5), R = D-2,\\
	&D=6, D-2=4, R\equiv 13~(\textup{mod} ~6), R \neq D-2.
\end{align*}
As seen here, in the ${2^{nd}}$ and ${4^{th}}$ case, the constraints on ${R}$ are not fulfilled.
Therefore, according to the discussion above, ${2}$ will be present in the slack list of ${13}$. Thus, ${13+2=15}$ cannot be a prime number according to the proposed method. It is well known that  ${15}$ is not a prime number.
\\
Going further, the mathematical proof for all twin primes can be easily derived from eq.(11). If ${R=D-2}$ for any division ${\dfrac{P}{D}}$, then slack, ${S=2}$. Thus, ${P+2}$ will not be prime, according to the proposed algorithm. Therefore, ${P}$ cannot have a twin prime as its successor. An added verification was done by writing a computer program to generate twin primes within different ranges by providing an upper bound, ${U}$, and it was found that the program successfully generates twin primes within ${U}$. The program detects twin primes on the basis of the \emph{{R-constraints}} mentioned earlier. 
Therefore, if any prime number has to be tested for presence of a twin as its successor prime, it has to fulfill each constraint corresponding to one of its potential divisors. If it does, it will turn out to be the first number of the twin prime pair, ${\left(P_{i},P_{i+1}\right)}$.
The following observations can be made at this point:
\begin{enumerate}
   \item Since there are infinitely many primes, and for any   two primes, ${\left(P_{1},P_{2}\right)}$, ${ m\equiv P_{1} ~(\textup{mod} ~D)\equiv P_{2} ~(\textup{mod} ~D)}$ may hold, the remainder ${R}$ can have the same values for different prime numbers ${\left(P_{1},P_{2},...,P_{n}\right)}$, corresponding to a given divisor ${D}$.
    \item As one successively moves through the infinite set of primes(${P}$), the number divisors(${D}$) goes on increasing. Every division, ${\dfrac{P}{D}}$ will result in a remainder, ${R}$. If all remainders obtained satisfy ${R \neq D-2}$, ${P}$ has a twin prime, else not. However, it cannot be predetermined which of the primes will simultaneously satisfy all the inequalities mentioned in ${R-constraints}$. 
\end{enumerate}

\section{Conclusion}
This paper presented a new method to calculate the next prime number after a given prime, ${P}$. This method helped to establish a rigid upper bound on the maximal gap ${G_{i}}$ between two consecutive primes,${\left(P_{i},P_{i+1}\right)}$ as 
\begin{equation*}
G_{i}\leq\lfloor\dfrac{P_{i}+1}{2}\rfloor,    
\end{equation*}
which is much smaller than that set by Bertrand's postulate, ${G_{i}<P_{i}}$(eq. 2).
Also, a new property of twin primes was unfolded in the form of the ${R-constraints}$, which need to be satisfied by the first number of any twin prime pair. 
In order to prove the \textbf{Twin Prime Conjecture}, 
\begin{equation*}
\lim_{n \to \infty}\left(P_{n+1}-P_{n}\right)=2,    
\end{equation*}
it would need to be proven that this system of constraints/inequalities has infinite solutions.
Therefore, this paper makes some contributions towards the understanding of prime numbers, prime gaps, twin primes and the very notable Twin Prime Conjecture.






\bibliographystyle{plain}
\bibliography{sample.bib}







\end{document}